\documentclass{amsart}
\theoremstyle{definition}

\theoremstyle{remark}

\numberwithin{equation}{section}

\def\A{{\mathcal A}}

\begin{document}

\title{On the complexity of Boolean matrix ranks}

\author{Yaroslav Shitov}
\address{National Research University Higher School of Economics, 20 Myasnitskaya Ulitsa, Moscow 101000, Russia}
\email{yaroslav-shitov@yandex.ru}


\begin{abstract}
We construct a reduction which proves that
the fooling set number and the determinantal rank
of a Boolean matrix are NP-hard to compute.
\end{abstract}

\maketitle

This note is devoted to the functions of \textit{determinantal rank} and \textit{fooling set number},
which are receiving attention in different applications, see~\cite{dSc, Theis} and references therein.
The purpose of this note is to give an NP-completeness proof for those functions, thus answering
questions asked in~\cite[Chapter 6]{dSc} and~\cite{Theis}.

A Boolean matrix $A\in\{0,1\}^{n\times n}$ is called \textit{sign-nonsingular} if
the condition $A_{1,\tau(1)}$$=$$\ldots$$=$$A_{n,\tau(n)}$$=$$1$ holds for some permutation $\tau$ on $\{1,\dots,n\}$ and for no permutation
of parity different from that of $\tau$. The \textit{determinantal rank} of a matrix is the size of its largest sign-nonsingular square submatrix.
Entries $(i_1,j_1),\dots,(i_n,j_n)$ of a Boolean matrix $B$ are called a \textit{fooling set} if $B_{i_s,j_t}=B_{i_t,j_s}=1$
holds exactly when $t=s$. The \textit{fooling set number}, equal to the cardinality of the largest fooling set,
is clearly an upper bound for the determinantal rank.

Consider a directed graph $G=(V,E)$ in which $(u,v)\in E$ implies $(v,u)\in E$ and construct the matrix $\A=\A(G)$ with rows
and columns indexed by $V\cup E$ as follows. Given $x,y\in V\cup E$, we set $\A_{xy}=1$ when $x=y$ or $(x,y)\in E$, or else
$y=(x,u)$ or $x=(u,y)$, for some vertex $u$; otherwise, we set $\A_{xy}=0$. If $U$ is an independent set in $G$, then by our
construction the submatrix of $A$ with rows and columns indexed by $U\cup E$ is sign-nonsingular. On the other hand, consider a subset
$S\subset (V\cup E)\times(V\cup E)$ such that $|S|=|E|+h$. If $A_{ij}$ with $(i,j)$ running over $S$ is a fooling set,
then the intersection of $S$ with $\{e,u\}\times\{e,v\}$ has cardinality $0$ or $1$, for every $e=(u,v)\in E$. So $S$ contains
at least $h$ elements of the form $(w,w)$ with $w\in V$, and those vertices $w$ form an independent set of cardinality $h$ in $G$.
The fooling set number and determinantal rank of $\A(G)$ are both therefore equal to $|U'|+|E|$, where $U'$ is the largest
independent set in $G$.

Checking whether a graph has an independent set of size exceeding $k$ is a classical NP-complete problem,
so it is NP-hard to decide whether or not the determinantal rank or fooling set number of a matrix
exceed $k$. Those problems are in fact NP-complete --- a polynomial algorithm for sign-singularity is a deep
result~\cite{RST} while yet a straightforward algorithm recognizing fooling sets is polynomial.

\end{document}